\documentclass[12pt]{article}
\usepackage{graphicx}
\usepackage{amsmath,amsthm,latexsym,amsbsy}
\newtheorem{definition}{Definition}
\newtheorem{theorem}{Theorem}
\newtheorem{remark}{Remark}

\newtheorem{lemma}{Lemma}
\newcommand{\h}{\hspace*{.24in}}
\def\geqslant {\ge}
\def\leqslant {\le}
\def\bq{\begin{equation}}
\def\eq{\end{equation}}
\title{Determine the source term of a two-dimensional heat equation}
\author{\normalsize DANG DUC TRONG $^1$, ~TRUONG TRUNG
TUYEN $^2$, \\\normalsize PHAN THANH NAM $^1$ and ALAIN PHAM NGOC DINH $^3$\\
\\\small\it $^1$ Mathematics Department, Natural Science HoChiMinh City University, Viet Nam
\\\small\it $^2$ Mathematics Department, Indiana University, Rawles Hall , Bloomington, IN 47405
\\\small\it $^3$ Mathematics Department, Mapmo UMR 6628, BP 67-59, 45067 Orleans cedex, France}
\date{{}}
\begin{document}
\maketitle
\begin{abstract} Let $\Omega$ be a two-dimensional heat conduction body. We consider the problem of determining the heat source $F(x,t)=\varphi(t)f(x,y)$ with $\varphi$ be given inexactly and $f$ be unknown. The problem is nonlinear and ill-posed. By a specific form of Fourier transforms, we shall show that the heat source is determined uniquely by the minimum boundary condition and the temperature distribution in $\Omega$ at the initial time $t=0$ and at the final time $t=1$. Using the methods of Tikhonov's regularization and truncated integration, we construct the regularized solutions.
\\\h MSC 2000: 35K05, 42B10, 65M32.
\\\h Key words: Error estimate, Fourier transform, ill-posed problem, heat source, Tikhonov's regularization, truncated integration.
\end{abstract}
\text{}\\
{\bf 1. Introduction and main results}\\\\
\h Let $\Omega=(0,1)\times (0,1)$ be a heat conduction body and $u(x,y,t)$ be the temperature in $\Omega$. We consider the problem of determining a pair of functions $(u,f)$ satisfying the system
\bq
\left\{ \begin{gathered}
  u_t  - \Delta u = \varphi (t)f(x,y),\h t\in (0,1); (x,y)\in \Omega, \hfill \\
  u_x (0,y,t) = u_x (1,y,t) = u_y (x,0,t) = u_y (x,1,t) = 0, \hfill \\
  u(1,y,t) = 0 ,\hfill \\
 \end{gathered}  \right.\label{1}
\eq
subject to the initial datum and the final datum
$$u(x,y,0) = g_0,u(x,y,1) =g_1.$$
Here, $\varphi\in L^1(0,T)$ and $g_0,g_1\in L^2(\Omega)$ are given inexactly.
\\\h This is a case of the problem of finding heat source $F(\xi,t,u)$ satisfying the heat equation
\[
u_t  - \Delta u = F,
\]
for $\xi$ is the spacial variable. The problem has been investigated intensively for the last three decades by many authors. Because the problem is severely ill-posed and difficult, many preassumptions on the form of the heat source are required. Recently, in  \cite{CE1,CE2,Y1,STY,CY}, the authors reduced the heat source $F$ to the function that has separated form
$$F(\xi,t,u)=\varphi(t)f(\xi)$$
with $\varphi$ is unknown. Then, the authors in \cite{TLD, TQD} studied the problem in the case that the heat source $F(\xi,t,u)=\varphi(t)f(\xi)$ with $f$ is unknown. From some assumptions of $\varphi$, the author used the Fourier transform and truncated integration to regularize the problem with nonsmooth data. However, in \cite{TQD}, the Cauchy datum $u(x,y,t)$ is given on two parts of the boundary, i.e  the boundary $x=1$ and the boundary $y=1$, say. In the present paper, the Cauchy datum is given only on the boundary $x=1$. The requirement of the Cauchy datum is minimum because if the condition $u(1,y,t)$ is omitted then the uniqueness of the solution of the problem cannot hold. For example, we consider the system
$$
\left\{ \begin{gathered}
  u_t  - \Delta u = \varphi (t)f(x,y), \hfill \\
  u_x (0,y,t) = u_x (1,y,t) = u_y (x,0,t) = u_y (x,1,t) = 0, \hfill \\
\end{gathered}  \right.\label{ptt}
$$
subject to
$$  u(x,y,0) = g_0,u(x,y,1) = g_1. $$
This system has not the uniqueness property. Indeed, if
$$g_0=g_1=0, \varphi(t)=\pi cos(\pi t)+2\pi^2 \sin(\pi t)$$
then this system has, beside the trivial solution $(u,f)=(0,0)$, a nontrivial solution
\[
u(x,y,t) = \sin (\pi t)\cos (\pi x)\cos (\pi y),f(x,y) = \cos (\pi x)\cos (\pi y).
\]
\h In fact, by a specific form of Fourier transforms, we shall get
\begin{theorem} (Uniqueness) Let $u_1,u_2\in C^1([0,1];L^2(\Omega))\cap L^2(0,1;H^2(\Omega))$, $f_1,f_2\in L^2(\Omega)$. If $(u_i,f_i)$ ($i=1,2$) satisfy the system $(1)$ with $g_0,g_1\in L^2(\Omega)$  and $\varphi \in L^1 (0,1)\backslash \{ 0\}$ then
$$
(u_1 ,f_1 ) = (u_2 ,f_2 ).
$$
\end{theorem}
We also have a regularization result. Using the Tikhonov regularization and truncated integration, we can construct a regularized solution for all $\varphi \not \equiv 0$.
\begin{theorem} (Regularization 1) Let $(u_{ex},f_{ex})\in (C^1([0,1];L^2(\Omega))\cap L^2(0,1;H^2(\Omega)), L^2(\Omega))$ be the exact solution of the system $(1)$ corresponding the exact data $g_{0ex},g_{1ex}\in L^2(\Omega)$ and $\varphi_{ex}\in L^1 (0,1)\backslash \{ 0\}$. Let $g_{0\varepsilon},g_{1\varepsilon}\in L^2(\Omega)$ and $\varphi _{\varepsilon}\in L^1 (0,1)$ be measured data satisfying
\[
\left\| {g_{0\varepsilon}   - g_{0ex} } \right\|_{L^2 (\Omega )}  \leqslant \varepsilon ,\left\| {g_{1\varepsilon}   - g_{1ex} } \right\|_{L^2 (\Omega )}  \leqslant \varepsilon ,\left\| {\varphi _\varepsilon   - \varphi _{ex} } \right\|_{L^1 (0,1)}  \leqslant \varepsilon .
\]
\h From $\{g_{0\varepsilon},g_{1\varepsilon},\varphi_{\varepsilon}\}$, we can construct a regularized solution  $f_{1\varepsilon}\in L^2(R^2)$ such that
$$\mathop {\lim }\limits_{\varepsilon  \to 0} \left\| {f_{1\varepsilon}   - f_{ex} } \right\|_{L^2 (\Omega )}=0.$$
\h Moreover, if $f_{ex}\in H^1(\Omega)$ then for each $\beta\in (0,1/2)$, there exists $\varepsilon_0>0$ (depended on $u_{ex}, \varphi_{ex}$ and $\beta$) such that
\[
\left\| {f_{1\varepsilon}   - f_{ex} } \right\|^2_{L^2 (\Omega )} \le  4\pi ^{ - 1} \left( {1 + 9\left\| {f_{ex} } \right\|_{H^1 (\Omega )}^2 } \right)(ln(\varepsilon^{-1}))^{-\beta}
\]
for $0<\varepsilon<\varepsilon_0$.
\end{theorem}
A smoother regularized solution will be given if $\varphi$ satisfies the following condition
\vspace{0.1in}\text{}\\
$(H)$ There exist $\lambda\in(0,1)$ and $C_0>0$ such that either $\varphi(t)\ge C_0$ for a.e $t\in (\lambda,1)$ or $\varphi(t)\le -C_0$ for a.e $t\in (\lambda,1)$.
\vspace{0.1in}\text{}\\
\h We note that $(H)$ will be satisfied if $\varphi$ is continuous at $t=1$ and $\varphi(1)\ne 0$. Under $(H)$, one has
\begin{theorem} (Regularization 2) Let $(u_{ex},f_{ex})\in (C^1([0,1];L^2(\Omega))\cap L^2(0,1;H^2(\Omega)), L^2(\Omega))$ be the exact solution of the system $(1)$ corresponding the exact data $g_{0ex},g_{1ex}\in L^2(\Omega)$ and $\varphi_{ex}$ satisfying $(H)$. Let $g_{\varepsilon}=(g_{0\varepsilon},g_{1\varepsilon})\in (L^2(\Omega))^2$ and $\varphi _{\varepsilon}\in L^1 (0,1)$ be measured data satisfying
\[
\left\| {g_{0\varepsilon}   - g_{0ex} } \right\|_{L^2 (\Omega )}  \leqslant \varepsilon ,\left\| {g_{1\varepsilon}   - g_{1ex} } \right\|_{L^2 (\Omega )}  \leqslant \varepsilon ,\left\| {\varphi _\varepsilon   - \varphi _{ex} } \right\|_{L^1 (0,1)}  \leqslant \varepsilon .
\]
\h  From $\{g_{0\varepsilon},g_{1\varepsilon},\varphi_{\varepsilon}\}$, we can construct a regularized solution  $f_{2\varepsilon}\in L^2(R^2)$ such that
$$\mathop {\lim }\limits_{\varepsilon  \to 0} \left\| {f_{2\varepsilon}   - f_{ex} } \right\|_{L^2 (\Omega )}=0.$$
\h Moreover, if $f_{ex}\in H^1(\Omega)$ then there exist $\gamma>0$ (depended on $\varphi_{ex}$) and  $\varepsilon_0>0$ (depended on $u_{ex}, \varphi_{ex}$) such that
\[
\left\| {f_{2\varepsilon}   - f_{ex} } \right\|_{L^2 (\Omega )}^2  \leqslant 4\pi ^{ - 1} \left( {\left\| {f_{ex} } \right\|_{L^1 (\Omega )}^2 \varepsilon ^\gamma   + 9\left\| {f_{ex} } \right\|_{H^1 (\Omega )}^2 \varepsilon ^{1/6} } \right)
\]
for $0<\varepsilon<\varepsilon_0$.
\end{theorem}
The remainder of the paper is divided into three sections. In Section 2, we shall give some notations and preparation results. The main results will be  proven in Section 3. In Section 4, a numerical experiment will be given to illustrate our approximation.
\text{}\\\\
\text{\bf 2. Notations and preparation results}\\\\
\h First, we have
\begin{lemma} If $u\in C^1([0,1];L^2(\Omega))\cap L^2(0,1;H^2(\Omega))$ and $f\in L^2(\Omega)$ satisfy the system $(1)$ then for all $(\alpha,n)\in R\times Z$ we have
\[
\begin{gathered}
  \int\limits_\Omega  {\left( {g(x,y)-e^{-(\alpha ^2  + n^2 \pi ^2) }g(x,y)} \right)\cos (\alpha x)\cos (n\pi y)dxdy}\hfill\\
= \int\limits_0^1 {e^{(\alpha ^2  + n^2 \pi ^2 )(t-1)} \varphi (t)dt} .\int\limits_\Omega  {f(x,y)\cos (\alpha x)\cos (n\pi y)dxdy}.  \hfill \\
\end{gathered}
\]
\end{lemma}
\begin{proof} Getting the inner product (in $L^2(\Omega)$) the first equation of the system $(1)$ with $W=cos(\alpha x)cos(n\pi y)$, we have
\[
\frac{d}
{{dt}}\int\limits_\Omega  {uWdxdy}  + (\alpha ^2  + n^2 \pi ^2)\int\limits_\Omega  {uWdxdy} = \varphi (t)\int\limits_\Omega  {fWdxdy} .
\]
\h We multiply the latter equation with $e^{(\alpha ^2  + n^2 \pi ^2)(t-1)}$ to get
\bq
\begin{gathered}
   \frac{d}
{{dt}}\left({e^{(\alpha ^2  + n^2 \pi ^2 )(t-1)}  \int\limits_\Omega  {uWdxdy} } \right) = e^{(\alpha ^2  + n^2 \pi ^2 )(t-1)} \varphi (t)\int\limits_\Omega  {fWdxdy} \label{tam}.\hfill \\
\end{gathered}
\eq
\h Integrating (from 0 to 1) the latter equality with respect to $t$, we shall get the desired result.
\end{proof}
\begin{definition} For all $\varphi  \in L^1 (0,1)$, we set $D(\varphi):R\times Z\to R$
$$
D(\varphi )(\alpha ,n) = \int\limits_0^1 {e^{(\alpha ^2  +n^2\pi ^2 )(t-1)} \varphi (t)dt},
$$
and
\[
B(\varphi ,r,\sigma) = \left\{ {\alpha  \in ( - r,r)\left| {\exists n \in Z \cap ( - r,r),\left| {D(\varphi )(\alpha ,n)} \right| \leqslant \sigma} \right.} \right\}.
\]
\end{definition}
Let $A$ be a subset of $R$. From now on we denote by $m(A)$ the Lebesgue measure of $A$. Using the idea in \cite{TT} (see Theorem 4), we have the following result.
\begin{lemma} Let $\varphi  \in L^1 (0,1)$, $\beta \in (0,1/2)$ and $q>q_1>0$.
\\$(i)$ Assume that $\varphi \not \equiv 0$. Then $D(\varphi )(\alpha,n)\ne 0$ for all $n\in Z$ and for a.e $\alpha \in R$. Moreover, there exists $\varepsilon_0>0$ (depended on $\varphi, q$ and $\beta$) such that
$$m(B(\varphi ,R_{1\varepsilon},\varepsilon^q ))<R^{-1}_{1\varepsilon}$$
for $0<\varepsilon<\varepsilon_0$ and $R_{1\varepsilon }  := (\ln (\varepsilon ^{ - 1}))^{\beta}$.
\\$(ii)$ Assume that $\varphi$ satisfies (H). Then there exist $\gamma>0$ (depended on $\varphi$) and $\varepsilon_0>0$ (depended on $\varphi, q$ and $q_1$) such that
$$m(B(\varphi ,R_{2\varepsilon},\varepsilon^q ))<\varepsilon^{\gamma}$$
for $0<\varepsilon<\varepsilon_0$ and $R_{2\varepsilon}:= \varepsilon^{-q_1/2}$.
\end{lemma}
\begin{proof} Because $\varphi  \in L^1 (0,1)\backslash \{ 0\}$, the map $\phi:C\to C$
\[
\phi(z)=e^{-z}\int\limits_0^1 {e^{zt} \varphi (t)dt}=
\int\limits_0^1 {e^{-zt} \varphi (1-t)dt}
\]
is a nontrivial entire function. Hence, for each $n\in Z$, the function
$$
\phi_n(z)=\phi(z^2+n^2\pi^2)=\int\limits_0^1 {e^{(z^2+n^2\pi^2)(t-1)} \varphi (t)dt}
$$
is also a nontrivial entire function.
\\$(i)$ For each $n\in Z$, since the zeros set of $\phi_n$ is either finite or countable, $D(\varphi )(\alpha ,n)=\phi_n(\alpha)\ne 0$ for a.e $\alpha\in R$. Hence $D(\varphi )(\alpha,n)\ne 0$ for all $n\in Z$ and for a.e $\alpha \in R$.
\\\h To estimate the measure of the set $B(\varphi ,r,\sigma)$, we shall use the following result in \cite{Ya} (Theorem 4 of Section $11.3$).
\begin{lemma} Let $f(z)$ be a function analytic in the disk $\{z:|z|\le 2eR\}$, $|f(0)|=1$, and let $\eta$ be an arbitrary small positive number. Then the estimate
\[
\ln |f(z)| >  - \ln (\frac{{15e^3 }}
{\eta }).\ln (M_f (2eR))
\]
is valid everywhere in the disk $\{z:|z|\le R\}$ except a set of disks $(C_j)$ with sum of radii $\sum r_j\le \eta R$. Here $M_f (r) = \mathop {\max }\limits_{|z| = r} |f(z)|$.
\end{lemma}
Since $\phi\not\equiv 0$, there exists $a_0\in (-1,0)$ such that $|\phi(a_0)|=C_1>0$. For each $n\in Z$, $|n|\le R_{1\varepsilon}$, we put $z_n=i\sqrt{n^2\pi^2-a_0}$. Then $\left| {z_n } \right| \leqslant \pi R_{1\varepsilon }  + 1
$ and $|\phi_n(z_n)|=|\phi(a_0)|=C_1$. We have that
\[
\Psi_n (z) := \frac{{\phi_n (z + z_n )}} {{C_1}}
\]
is an entire function and $|\Psi_n(0)|=1$, moreover for all $z \in C,\left| z \right| \leqslant 2eR$, $n\in Z, |n|\le R_{1\varepsilon}$,
\[
\begin{gathered}
  \left| {\Psi _n (z)} \right| =\frac{1}{C_1} .\left| {\int\limits_0^1 {e^{\left( {(z + z_n )^2  + n^2\pi^2 } \right)(t-1)} \varphi (t)dt} } \right| \hfill\\
\leqslant \frac{1}{C_1} .\int\limits_0^1 {e^{(\left| z \right|+\left| z_n \right|)^2} \left| {\varphi (t)} \right|dt}
\le  e^{(2eR+R_{1\varepsilon}\pi+1)^2}.\frac{\left\| \varphi  \right\|_{L^1 (0,1)} }{C_1} .\hfill \\
\end{gathered}
\]
\h For $\varepsilon>0$ small enough (depended on $\varphi$, $q$ and $\beta$) and for each $n\in Z, |n|\le R_{1\varepsilon}$, applying Lemma 3 to $R = (1+\pi)R_{1\varepsilon}+1$ and $\eta  = \frac{1}{5RR^2_{1\varepsilon}}$, we obtain that
\[
\begin{gathered}
  \ln\left| {\Psi _n (z)} \right| >  - \ln (\frac{{15e^3 }}
{\eta }).\ln \left( {M_{\Psi _n } (2eR)} \right) \hfill \\
   \geqslant  - \left[ {\ln (R) + 2\ln (R_{1\varepsilon } ) + \ln (75e^3 )} \right].\left[ {(2eR + \pi R_{1\varepsilon }  + 1)^2  + \ln \left( {\frac{{\left\| \varphi  \right\|_{L^1 (0,1)} }}
{{C_1 }}} \right)} \right] \hfill \\
   \geqslant  - 3\ln (R_{1\varepsilon } )(2e(1 + \pi ) + \pi  + 1)^2 R_\varepsilon ^2  + \ln \left( {\frac{1}
{{C_1 }}} \right) \hfill \\
   \geqslant  - q\ln \left( {\varepsilon ^{ - 1} } \right) + \ln \left( {\frac{1}
{{C_1 }}} \right) = \ln \left( {\frac{{\varepsilon ^q }}
{{C_1 }}} \right) \hfill \\
\end{gathered}
\]
for all $|z|\le R$  except a set of disks $\{B(z_{nj},r_{nj})\}_{j\in J_n}$ with sum of radii
$$\sum\limits_{j \in J_n } {r_{nj} }  \leqslant \eta R = \frac{1}
{{5R_{1\varepsilon }^2 }}.$$
\\\h Consequently, for $\varepsilon>0$ small enough and for each $n\in Z, |n|\le R_{1\varepsilon}$  we get
\[
\left| {D(\varphi )(\alpha ,n)} \right|=\left| {\phi_n (\alpha)} \right| = C_1 .\left| {\Psi _n (\alpha - z_n )} \right| > C_1 .\frac{1}
{{C_1 }}\varepsilon ^q  = \varepsilon ^q
\]
for all $|\alpha|  \in [-R_{1\varepsilon},R_{1\varepsilon}]$ except the set $\mathop  \cup \limits_{j \in J_n} B(z_{nj}  + z_n ,r_{nj} )$. So
$$B(\varphi ,R_{1\varepsilon } ,\varepsilon ^q ) \subset \mathop  \cup \limits_{n\in Z, |n|\le R_{1\varepsilon}}\mathop  \cup \limits_{j \in J_n} (\xi_{nj}-r_{nj},\xi_{nj}+r_{nj})$$
with $\xi_{nj}=Re(z_{nj})$. Thus
\[
m(B(\varphi ,R_{1\varepsilon } ,\varepsilon ^q )) \leqslant \sum\limits_{ - R_{1\varepsilon }  < n < R_{1\varepsilon } } {\sum\limits_{j \in J_n } {2r_{nj} } }  \leqslant (2R_{1\varepsilon }  + 1).2.\frac{1}
{{5R_{1\varepsilon }^2 }} < \frac{1}
{{R_{1\varepsilon } }}
\]
for $\varepsilon>0$ small enough (depended on $\varphi$, $q$ and $\beta$).
\\$(ii)$  Note that
\[
\begin{gathered}
  D(\varphi )(\alpha ,n) = \int\limits_0^\lambda  {e^{(\alpha ^2  + n^2 \pi ^2 )(t - 1)} \varphi (t)dt}  + \int\limits_\lambda ^1 {e^{(\alpha ^2  + n^2 \pi ^2 )(t - 1)} \varphi (t)dt}  \hfill \\
\h\h\h   \geqslant  - \int\limits_0^\lambda  {e^{(\alpha ^2  + n^2 \pi ^2 )(\lambda  - 1)} \left| {\varphi (t)} \right|dt}  + \int\limits_\lambda ^1 {e^{(\alpha ^2  + n^2 \pi ^2 )(t - 1)} C_0 dt}  \hfill \\
\h\h\h   \geqslant  - e^{(\alpha ^2  + n^2 \pi ^2 )(\lambda  - 1)} .\left\| \varphi  \right\|_{L^1 (0,1)}  + C_0 .\frac{{1 - e^{(\alpha ^2  + n^2 \pi ^2 )(\lambda  - 1)} }}
{{\alpha ^2  + n^2 \pi ^2 }}. \hfill \\
\end{gathered}
\]
Therefore, there exists a constant $R_1>0$ (depended on $\varphi, \lambda$) satisfying for either $|\alpha|\ge R_1$ or $|n|\ge R_1$ that
\[
D(\varphi )(\alpha ,n) \geqslant \frac{{C_0 }}
{{2(\alpha ^2  + n^2 \pi ^2 )}}.
\]
Consequently, for $\varepsilon>0$ small enough (depended on $\varphi$, $q$ and $q_1$) and for all $(\alpha,n)\in (-R_{2\varepsilon},R_{2\varepsilon})^2\backslash (-R_1,R_1)^2$, we get
$$
D(\varphi )(\alpha ,n) \geqslant\frac{{C_0 }}
{{2(\alpha ^2  + n^2 \pi ^2 )}} \ge  \frac{C_0}{2(1+\pi^2)R^2_{2\varepsilon}}=\frac{C_0}{2(1+\pi^2)}.\varepsilon^{q_1}>\varepsilon^{q}.
$$
\h Now we consider only $(\alpha,a)\in (-R_1,R_1)^2$. Let $a_0$, $z_n$, $\Psi_n$ as in $(i)$ and put $R_2=(1+\pi)R_1+1$. Then for all $z \in C,\left| z \right| \leqslant 2eR_2$, we get
\[
\begin{gathered}
  \left| {\Psi _n (z)} \right| =\frac{1}{C_1} .\left| {\int\limits_0^1 {e^{\left( {(z + z_n )^2  + n^2\pi^2 } \right)(t-1)} \varphi (t)dt} } \right|
\le  e^{(2eR_2+R_{1}\pi+1)^2}.\frac{\left\| \varphi  \right\|_{L^1 (0,1)} }{C_1}\le C_2 \hfill \\
\end{gathered}
\]
where $C_2>1$ be a constant independent of $n$.
\\\h For $\varepsilon>0$ small enough, applying Lemma 3 to $R = R_2$, $\eta  =\dfrac{ \varepsilon^{\gamma}}{(4R_1+2)R_2}$ and $\gamma=\dfrac{q_1}{2ln(C_2)}>0$, we get
\[
\begin{gathered}
  \ln \left| {\Psi_n (z)} \right| > \left[ {\gamma .\ln (\varepsilon ) - \ln (75R_2e^3)} \right].\ln (C_2 )
   > 2\gamma \ln (C_2 )\ln (\varepsilon )-\ln(C_1)= \ln (C_1^{-1}\varepsilon ^{q_1} ) \hfill \\
\end{gathered}
\]
for all $|z|\le R_2$ except a set of disks $\{B(z_{nk},r_{nk})\}_{k\in K_n}$ with sum of radii $$
\sum\limits_{k \in K_n } {r_{nk} }  \leqslant \eta R_2 = \frac{\varepsilon^\gamma}
{{4R_1+2 }}.$$
Consequently, for $\varepsilon>0$ small enough and for each $n\in Z, |n|\le R_{1}$, we have
\[
\left| {D(\varphi )(\alpha ,n)} \right|=\left| {\phi_n (\alpha)} \right| = C_1 .\left| {\Psi _n (\alpha - z_n )} \right| >\varepsilon ^{q_1}>\varepsilon ^q
\]
for all $\alpha\in [-R_{1},R_1]$ except the set $\mathop  \cup \limits_{k \in K_n} B(z_{nk}  + z_n ,r_{nk} )$. Thus
$$B(\varphi,R_{2\varepsilon},\varepsilon^q) \subset \mathop  \cup \limits_{n\in Z, |n|\le R_{1}}\mathop  \cup \limits_{k \in K_n} (\xi_{nk}-r_{nk},\xi_{nk}+r_{nk})$$
with $\xi_{nk}=Re(z_{nk})$. Therefore,
\[
m(B(\varphi ,R_{2\varepsilon } ,\varepsilon ^q )) \leqslant \sum\limits_{ - R_{1}  < n < R_{1} } {\sum\limits_{k \in K_n } {2r_{nk} } }< (2R_{1}  + 1).2.\frac{\varepsilon^\gamma}
{{4R_{1}+2}} =\varepsilon^\gamma
\]
for $\varepsilon>0$ small enough (depended on $\varphi$, $q$ and $\beta$).
\\\h The proof of Lemma 2 is completed.
\end{proof}
\begin{definition} For each $w\in L^2(\Omega)$, we set $G(w)$ defined on $R\times Z$ by
\[
G(w)(\alpha,n)=\int\limits_\Omega  {w(x,y)\cos (\alpha x)\cos (n\pi y)dxdy}.
\]
\end{definition}
\begin{lemma} $(i)$ For each $v\in L^2(0,1)$ we have
\[
\int\limits_{ - \infty }^\infty  {\left| {\int\limits_0^1 {v(x)\cos (\alpha x)dx} } \right|^2 d\alpha}  = \pi \left\| v \right\|_{L^2 (0,1)}^2 .
\]
$(ii)$ For each $w\in L^2(\Omega)$ we have
$$
 \sum\limits_{n =  - \infty }^\infty  {\int\limits_{ - \infty }^\infty  {\left| {{G(w)(\alpha,n)} } \right|^2 d\alpha } }  = \pi\left\| w \right\|_{L^2 (\Omega )}^2 .
$$
\end{lemma}
\begin{proof} $(i)$ Putting
\[
\widetilde v(x) =\left\{ \begin{gathered}
  v(x),\h x\in(0,1), \hfill \\
  v( - x),\h x\in (-1,0) ,\hfill \\
  0 ,\h x\notin (-1,1),\hfill \\
\end{gathered}  \right.
\]
then the Fourier transform of $\widetilde v(x)$ is
\[
F(\widetilde v)(\alpha ): = \int\limits_{ - \infty }^\infty  {\widetilde v(x)e^{ - i\alpha x} dx}  =2\int\limits_0^1 {v(x)\cos (\alpha x)}dx.
\]
Using Parseval equality, we get
\[
\int\limits_{ - \infty }^\infty  {\left| {\int\limits_0^1 {v(x)\cos (\alpha x)dx} } \right|^2 d\alpha}  =\frac{1}{4} \left\| {F(\widetilde v)} \right\|_{L^2 (R)}^2  = \frac{\pi}{2}\left\| {\widetilde v} \right\|_{L^2 (R)}^2  = \pi \left\| v \right\|_{L^2 (0,1)}^2 .
\]
$(ii)$  For each $n\in Z$, we put
\[
h_n (w)(x) = \int\limits_0^1 {w(x,y)\cos (n\pi y)dy} .
\]
Applying $(i)$ to $v=h_n(w)$, we get
\[
\begin{gathered}
  {\int\limits_{ - \infty }^\infty  {\left| { {G(w)(\alpha,n)} } \right|^2 d\alpha } }
   = \int\limits_{ - \infty }^\infty  {\left| {\int\limits_0^1 {h_n (w)(x)\cos (\alpha x)dx} } \right|^2 d\alpha }  = \pi \left\| {h_n (w)} \right\|_{L^2 (0,1)}^2 . \hfill \\
\end{gathered}
\]
On the other hand, since $w\in L^2(\Omega)$, for a.e $x\in (0,1)$ we have $w(x,.)\in L^2(0,1)$ and its Fourier-cosin series corresponding to variable $y$ is $h_n(w)(x)$. Using Parseval equality, we get
\[
\sum\limits_{n =  - \infty }^\infty  {\left| {h_n (w)(x)} \right|^2 }  = \left\| {w(x,.)} \right\|_{L^2 (0,1)}^2  = \int\limits_0^1 {\left| {w(x,y)} \right|^2 dy},\h a.e\h x\in (0,1) .
\]
Therefore,
\[
\begin{gathered}
  \sum\limits_{n =  - \infty }^\infty  {\int\limits_{ - \infty }^\infty  {\left| {{G(w)(\alpha,n)} } \right|^2 d\alpha } }  = \pi \sum\limits_{n =  - \infty }^\infty  {\left\| {h_n (w)} \right\|_{L^2 (0,1)}^2 }
   = \pi \sum\limits_{n =  - \infty }^\infty  {\left( {\int\limits_0^1 {\left| {h_n (w)(x)} \right|^2 dx} } \right)} \hfill \\
  = \pi \int\limits_0^1 {\left( {\sum\limits_{n =  - \infty }^\infty  {\left| {h_n (w)(x)} \right|^2 } } \right)dx}
   = \pi \int\limits_0^1 {\left( {\int\limits_0^1 {\left| {w(x,y)} \right|^2 dy} } \right)dx = \pi \left\| w \right\|_{L^2 (\Omega )}^2 }  .\hfill \\
\end{gathered}
\]
\h The proof is completed.
\end{proof}
To prove the regularization results, we need one more preparation.
\begin{definition} For each $w\in L^2(\Omega)$ and $r>0$, we set
\[
\mu (w ,r) = \sum\limits_{|n| \geqslant r  } {\int\limits_{ - \infty }^\infty  {\left| {G(w )(\alpha,n)} \right|^2 d\alpha }}+ \sum\limits_{n =  - \infty }^\infty  {\int\limits_{\left| \alpha  \right| \geqslant r } {\left| {G(w )(\alpha,n)} \right|^2 d\alpha } } .
\]
\end{definition}
\begin{lemma} For each $w\in L^2(\Omega)$, we have $\mathop {\lim }\limits_{r \to  + \infty } \mu (w,r) = 0$. Moreover, if $w\in H^1(\Omega)$ then
\[
\mu (w ,r ) \leqslant \left( {\frac{{8}} {{r  }} + \frac{2\pi } {{r
^2 }}} \right)\left\| {w} \right\|_{H^1 (\Omega )}^2 .
\]
\end{lemma}
\begin{proof} For each $w\in L^2(\Omega)$, applying Lemma 4 we obtain
$$
 \sum\limits_{n =  - \infty }^\infty  {\int\limits_{ - \infty }^\infty  {\left| {\int\limits_\Omega  {G(w)(\alpha,n)} } \right|^2 d\alpha } }  = \pi\left\| w \right\|_{L^2 (\Omega )}^2 <+\infty.
$$
It implies that $\mathop {\lim }\limits_{r \to  + \infty } \mu (w,r) = 0$.
\\\h Now, we consider $w\in H^1(\Omega)$. Since
\[
\begin{gathered}
  \int\limits_0^1 {w(x,y)\cos (n\pi y)dy}  \hfill \\
   = \left[ {w(x,y)\frac{{\sin (n\pi y)}}
{{n\pi }}} \right]_{y = 0}^{y = 1}  - \int\limits_0^1 {\frac{{\partial w}}
{{\partial y}}(x,y)\frac{{\sin (n\pi y)}}
{{n\pi }}dy}
   =  - \frac{1}
{{n\pi }}\int\limits_0^1 {\frac{{\partial w }}
{{\partial y}}(x,y)\sin (n\pi y)dy}  ,\hfill \\
\end{gathered}
\]
we get
\[
G(w)(\alpha,n) =  - \frac{1}
{{n\pi }}\int\limits_\Omega  {\frac{{\partial w }}
{{\partial y}}(x,y)\cos (\alpha x)\sin (n\pi y)dxdy} = - \frac{1}
{{n\pi }}G(\frac{{\partial w }}
{{\partial y}})(\alpha,n).
\]
Consequently,
\[
\begin{gathered}
  \sum\limits_{|n| \geqslant r  } {\int\limits_{ - \infty }^\infty  {\left| {G(w)} \right|^2 d\alpha } }
   \leqslant \frac{1}
{{r^2 \pi ^2 }}\sum\limits_{|n| \geqslant r  } {\int\limits_{ - \infty }^\infty  {\left| {G(\frac{{\partial w }}
{{\partial y}})(\alpha,n) } \right|^2 d\alpha } }  \hfill \\
\le \frac{1}
{{r^2 \pi ^2 }}\sum\limits_{n =  - \infty }^\infty {\int\limits_{ - \infty }^\infty  {\left| {G(\frac{{\partial w }}
{{\partial y}})(\alpha,n) } \right|^2 d\alpha } }
     = \frac{1}
{{r^2 \pi  }}\left\| {\frac{{\partial w }}
{{\partial y}}} \right\|_{L^2 (\Omega )}^2 .\hfill \\
\end{gathered}
\]
Similarly, we have
\[
\begin{gathered}
  \int\limits_0^1 {w (x,y)\cos (\alpha x)dx}  \hfill \\
   = \left[ {w(x,y)\frac{{\sin (\alpha x)}}
{\alpha }} \right]_{x = 0}^{x = 1}  - \int\limits_0^1 {\frac{{\partial w }}
{{\partial x}}(x,y)\frac{{\sin (\alpha x)}}
{\alpha }dx}  \hfill \\
   = \frac{{\sin (\alpha )}}
{\alpha }w (1,y) - \frac{1}
{\alpha }\int\limits_0^1 {\frac{{\partial w }}
{{\partial x}}(x,y)\sin (\alpha x)dx}.  \hfill \\
\end{gathered}
\]
Therefore,
\[
G(w )(\alpha,n) = \frac{{\sin (\alpha )}}
{\alpha }\int\limits_0^1 {w (1,y)\cos (n\pi y)dy}  - \frac{1}
{\alpha }G\left(\frac{{\partial w }}{{\partial x}}\right)(\alpha,n) .
\]
Hence
\[
\begin{gathered}
  \sum\limits_{n =  - \infty }^\infty  {\int\limits_{\left| \alpha  \right| \geqslant r } {\left| {G(w )} \right|^2 d\alpha } } \hfill\\
 \leqslant \int\limits_{\left| \alpha  \right| \geqslant r  } {\frac{2}
{{\alpha ^2 }}d\alpha } .\sum\limits_{n =  - \infty }^\infty
{\left(\int_0^1\left| {w (1,y)\cos (n\pi y)} \right|dy\right)^2 }
   + \frac{2}
{{r^2 }}\sum\limits_{n =  - \infty }^\infty  {\int\limits_{ - \infty }^\infty  {\left| {G\left(\frac{{\partial w }}{{\partial x}}\right) } \right|^2 d\alpha } }  \hfill \\
   = \frac{4}
{{r  }}\left\| {w (1,.)} \right\|_{L^2 (0,1)}^2  + \frac{2}
{{r^2 }}.\pi\left\| {\frac{{\partial w }}
{{\partial x}}} \right\|_{L^2 (\Omega )}^2.  \hfill \\
\end{gathered}
\]
Noting that
\[
w(1,y) = \int\limits_0^1 {\frac{\partial }
{{\partial x}}\left( {xw (x,y)} \right)dx}  = \int\limits_0^1 {\left( {w (x,y) + x\frac{{\partial w}}
{{\partial x}}(x,y)} \right)dx} ,
\]
we get
\[
\begin{gathered}
  \left| {w (1,y)} \right|^2  \leqslant \left| {\int\limits_0^1 {\left( {w (x,y) + x\frac{{\partial w }}
{{\partial x}}(x,y)} \right)dx} } \right|^2
   \leqslant 2\int\limits_0^1 {\left( {\left| {w (x,y)} \right|^2  + \left| {\frac{{\partial w }}
{{\partial x}}(x,y)} \right|^2 } \right)dx} , \hfill \\
\end{gathered}
\]
and
\[
\left\| {w (1,.)} \right\|_{L^2 (0,1)}^2  \leqslant 2\left\| {w} \right\|_{L^2 (\Omega )}^2  + 2\left\| {\frac{{\partial w }}
{{\partial x}}} \right\|_{L^2 (\Omega )}^2  \leqslant 2\left\| {w } \right\|_{H^1 (\Omega )}^2 .
\]
Thus
\[
\sum\limits_{n =  - \infty }^\infty  {\int\limits_{\left| \alpha  \right| \geqslant r  } {\left| {G(w )} \right|^2 d\alpha } }  \leqslant \frac{{8 }}
{{r  }}\left\| {w } \right\|_{H^1 (\Omega )}^2  + \frac{2\pi}
{{r ^2 }}\left\| {\frac{{\partial w }}
{{\partial x}}} \right\|_{L^2 (\Omega )}^2
\]
In summary, we get
\[
\mu (w ,r ) \leqslant \left( {\frac{{8}} {{r  }} + \frac{2\pi } {{r
^2 }}} \right)\left\| w \right\|_{H^1 (\Omega )}^2 .
\]
\h The proof is completed.
\end{proof}
\text{}\\
\text{\bf 3. Proofs of main theorems}\\
\text{}\\
{\bf Proof of Theorem 1}
\begin{proof} Put $u=u_1-u_2$, $f=f_1-f_2$. Then $(u,f)$ satisfies system $(1)$ corresponding to $g_0=g_1=0$. Let D and G be as in Definition 1 and Definition 2. Applying Lemma 1, we obtain
$$
D( \varphi)(\alpha,n).G(f)(\alpha,n)=0,\h\forall (\alpha,n) \in R\times Z.
$$
\h For each $n\in Z$, according to Lemma 2, $D( \varphi)(\alpha,n)\ne 0$ for a.e $\alpha\in R$. It implies $G(f)(\alpha,n)=0$ for a.e $\alpha\in R$. So $f\equiv 0$ because of Lemma 4. Hence, equation $(\ref{tam})$ in the proof of Lemma 1 becomes
\[
\frac{d}
{{dt}}\left( {e^{(\alpha ^2  + n^2 \pi ^2 )(t-1)} \int\limits_\Omega  {u(x,y,t)\cos (\alpha x)\cos (n\pi y)dxdy} } \right) = 0.
\]
Since $u(x,y,0)=0$, the latter equation implies that for all $t\in (0,1)$
\[
\int\limits_\Omega  {u(x,y,t)\cos (\alpha x)\cos (n\pi y)dxdy}  = 0
\]
Using Lemma 4 again, we obtain $u\equiv 0$ as desired.
\end{proof}
\text{}\\
{\bf Proof of Theorem 2 and Theorem 3}
\begin{proof} The proof is divided into four steps.
\text{}\\{\bf Step 1.} Let $D$ and $G$ be as in Definition 1 and Definition 2. For each $g=(g_0,g_1)\in (L^2(\Omega))^2$, we put
\[
H(g)(\alpha ,n) = G(g_1)(\alpha,n) - e^{-(\alpha ^2  + a^2 \pi ^2)} G(g_0 )(\alpha,n).
\]
From Lemma 1 and Lemma 2, it follows that
\[
G(f_{ex})(\alpha ,n) = \frac{{H(g_{ex} )(\alpha ,n)}}
{{D(\varphi_{ex})(\alpha ,n)}}
\]
for a.e $\alpha\in R$ and for all $n\in Z$. From Lemma 4, we have
\[
f_{ex} (x,y) = \frac{1}
{\pi }.\sum\limits_{n =  - \infty }^\infty  {\left( {\int\limits_{ - \infty }^\infty  {\frac{{H(g_{ex} )(\alpha ,n)}}
{{D(\varphi _{ex} )(\alpha ,n)}}\cos (\alpha x)d\alpha } } \right)\cos (n\pi y)} .
\]
\h We shall construct the regularized solution $f_{\varepsilon}$ by the following formula
\[
f_\varepsilon  (x,y) = \frac{1}
{\pi }.\sum\limits_{ - R_\varepsilon   < n < R_\varepsilon  } {\left( {\int\limits_{ - R_\varepsilon  }^{R_\varepsilon  } {\frac{{H(g_\varepsilon  )D(\varphi _\varepsilon  )}}
{{D^2 (\varphi _\varepsilon  ) + \delta _\varepsilon  }}\cos (\alpha x)d\alpha } } \right)\cos (n\pi y)} ,
\]
where $R_\varepsilon$ and $\delta _\varepsilon$ are regularized parameters chosen later.
\\\h It is obvious that $f_{\varepsilon}\in C(R^2)$ and
\[
G(f_\varepsilon  )(\alpha ,n) = \chi \left( {(( - R_\varepsilon  ,R_\varepsilon  )^2 } \right).\frac{{H(g_\varepsilon  )D(\varphi _\varepsilon  )}}
{{D^2 (\varphi _\varepsilon  ) + \delta _\varepsilon  }}
\]
where $\chi(A)$ is the characteristic function of the set $A$, i.e,
\[
\chi (A)(x) = \left\{ \begin{gathered}
  1,\h x \in A, \hfill \\
  0,\h x \notin A. \hfill \\
\end{gathered}  \right.
\]
\h To prove that $f_{\varepsilon}$ approximates $f_{ex}$ in $L^2(\Omega)$, we only have to verify that $G(f_{\varepsilon})$ approximates $G(f_{ex})$ in $L^2(R^2)$.
\\{\bf Step 2.} Estimate $\sum\limits_{n =  - \infty }^\infty  {\int\limits_{ - \infty }^\infty  {\left| {G(f_\varepsilon  ) - G(f_{ex} )} \right|^2 d\alpha } }
$.
\\\h We first estimate $|G(f_{\varepsilon})-G(f_{ex})|$. If $(\alpha,n)\notin (-R_\varepsilon,R_\varepsilon)^2$ then $\left| {G(f_\varepsilon  ) - G(f_{ex} )} \right| = \left| {G(f_{ex} )} \right| $. Now, we consider $(\alpha,n)\in (-R_\varepsilon,R_\varepsilon)^2$. Putting
\[
C_3  =\max \left\{ {\left\| {\varphi _{ex} } \right\|_{L^1 (0,1)}
,\left\| {g_{0ex} } \right\|_{L^2(\Omega )} ,\left\| {g_{1ex} }
\right\|_{L^2 (\Omega )} } \right\},\]
 by a direct calculation we
obtain
\[
\begin{gathered}
  \left| {H(g_{ex} )} \right| \leqslant \left\| {g_{1ex} } \right\|_{L^1 (\Omega )}  + \left\| {g_{0ex} } \right\|_{L^1 (\Omega )}  \leqslant 2C_3,\hfill \\
  \left| {H(g_\varepsilon  ) - H(g_{ex} )} \right| \leqslant \left\| {g_{1\varepsilon }  - g_{1ex} } \right\|_{L^1 (\Omega )}  + \left\| {g_{0\varepsilon }  - g_{0ex} } \right\|_{L^1 (\Omega )}  \leqslant 2\varepsilon , \hfill \\
  \left| {D(\varphi _{ex} )} \right| \leqslant \left\| {\varphi _{ex} } \right\|_{L^1 (0,1)}\le C_3, \left| {D(\varphi _\varepsilon  ) - D(\varphi _{ex} )} \right| \leqslant \left\| {\varphi _\varepsilon   - \varphi _{ex} } \right\|_{L^1 (0,1)}  \leqslant \varepsilon . \hfill \\
\end{gathered}
\]
We have
\[
\begin{gathered}
 \left| {G(f_\varepsilon  ) - G(f_{ex} )} \right| = \left| {\frac{{H(g_\varepsilon  )D(\varphi _\varepsilon  )}}
{{D^2 (\varphi _\varepsilon  ) + \delta _\varepsilon  }} - \frac{{H(g_{ex} )}}
{{D(\varphi _{ex} )}}} \right| \hfill \\
   \leqslant \left| {\frac{{H(g_\varepsilon  )D(\varphi _\varepsilon  )}}
{{D^2 (\varphi _\varepsilon  ) + \delta _\varepsilon  }} - \frac{{H(g_{ex} )D(\varphi _{ex} )}}
{{D^2 (\varphi _{ex} ) + \delta _{ex} }}} \right| + \left| {\frac{{H(g_{ex} )D(\varphi _{ex} )}}
{{D^2 (\varphi _{ex} ) + \delta _{ex} }} - \frac{{H(g_{ex} )}}
{{D(\varphi _{ex} )}}} \right| .\hfill \\
\end{gathered}
\]
We shall estimate each term of the right-hand side. We have
\[
\begin{gathered}
  \left| {\frac{{H(g_\varepsilon  )D(\varphi _\varepsilon  )}}
{{D^2 (\varphi _\varepsilon  ) + \delta _\varepsilon  }} - \frac{{H(g_{ex} )D(\varphi _{ex} )}}
{{D^2 (\varphi _{ex} ) + \delta _{ex} }}} \right| \hfill \\
   \leqslant \frac{{\left| {D(\varphi _\varepsilon  )D(\varphi _{ex} )} \right|.\left| {D(\varphi _{ex} )H(g_\varepsilon  ) - D(\varphi _\varepsilon  )H(g_{ex} )} \right|}}
{{\left( {D^2 (\varphi _\varepsilon  ) + \delta _\varepsilon  } \right)\left( {D^2 (\varphi _{ex} ) + \delta _{ex} } \right)}} + \frac{{\delta _\varepsilon  \left| {D(\varphi _\varepsilon  )H(g_\varepsilon  ) - D(\varphi _{ex} )H(g_{ex} )} \right|}}
{{\left( {D^2 (\varphi _\varepsilon  ) + \delta _\varepsilon  } \right)\left( {D^2 (\varphi _{ex} ) + \delta _{ex} } \right)}} \hfill \\
   \leqslant \frac{{\left| {D(\varphi _{ex} )H(g_\varepsilon  ) - D(\varphi _\varepsilon  )H(g_{ex} )} \right|}}
{{\delta _\varepsilon  }} + \frac{{\left| {D(\varphi _\varepsilon  )H(g_\varepsilon  ) - D(\varphi _{ex} )H(g_{ex} )} \right|}}
{{\delta _\varepsilon  }}. \hfill \\
\end{gathered}
\]
We get
\[
\begin{gathered}
  \left| {D(\varphi _{ex} )H(g_\varepsilon  ) - D(\varphi _\varepsilon  )H(g_{ex} )} \right| \hfill \\
   \leqslant \left| {D(\varphi _{ex} )} \right|.\left| {H(g_\varepsilon  ) - H(g_{ex} )} \right| + \left| {D(\varphi _\varepsilon  ) - D(\varphi _{ex} )} \right|.\left| {H(g_{ex} )} \right| \hfill \\
   \leqslant C_3.\left| {H(g_\varepsilon -g_{ex} )} \right|  +  \varepsilon.\left| {H(g_{ex} )} \right|, \hfill \\
\end{gathered}
\]
and similarly,
\[
\begin{gathered}
  \left| {D(\varphi _\varepsilon  )H(g_\varepsilon  ) - D(\varphi _{ex} )H(g_{ex} )} \right| \hfill \\
   \leqslant \left| {D(\varphi _{ex} )} \right|.\left| {H(g_\varepsilon  ) - H(g_{ex} )} \right| + \left| {D(\varphi _\varepsilon  ) - D(\varphi _{ex} )} \right|.\left| {H(g_{ex} )} \right| \hfill \\
  ~~ + \left| {\left( {D(\varphi _\varepsilon  ) - D(\varphi _{ex} )} \right)\left( {H(g_\varepsilon  ) - H(g_{ex} )} \right)} \right| \hfill \\
   \leqslant  2C_3.\left| {H(g_\varepsilon -g_{ex} )} \right|  +  \varepsilon.\left| {H(g_{ex} )} \right|.\hfill \\
\end{gathered}
\]
Thus,
\[
\left| {\frac{{H(g_\varepsilon  )D(\varphi _\varepsilon  )}}
{{D^2 (\varphi _\varepsilon  ) + \delta _\varepsilon  }} - \frac{{H(g_{ex} )D(\varphi _{ex} )}}
{{D^2 (\varphi _{ex} ) + \delta _\varepsilon  }}} \right| \le 3C_3.\left| {H(g_\varepsilon -g_{ex} )} \right|\delta _\varepsilon^{-1}  +  2\varepsilon.\left| {H(g_{ex} )} \right|\delta _\varepsilon^{-1}.
\]
We next consider
\[
\begin{gathered}
  \left| {\frac{{H(g_{ex} )D(\varphi _{ex} )}}
{{D^2 (\varphi _{ex} ) + \delta _\varepsilon  }} - \frac{{H(g_{ex} )}}
{{D(\varphi _{ex} )}}} \right| = \frac{{\delta _\varepsilon  }}
{{D^2 (\varphi _{ex} ) + \delta _\varepsilon  }}.\left| {\frac{{H(g_{ex} )}}
{{D(\varphi _{ex} )}}} \right|    = \frac{{\delta _\varepsilon  }}
{{D^2 (\varphi _{ex} ) + \delta _\varepsilon  }}.\left| {G(f_{ex} )} \right|. \hfill \\
\end{gathered}
\]
If $\left| {D(\varphi _{ex} )} \right| \geqslant \varepsilon ^q$ then
\[
\left| {\frac{{H(g_{ex} )D(\varphi _{ex} )}}
{{D^2 (\varphi _{ex} ) + \delta _\varepsilon  }} - \frac{{H(g_{ex} )}}
{{D(\varphi _{ex} )}}} \right| \leqslant \frac{{\delta _\varepsilon  }}
{{\varepsilon ^{2q} }}.\left| {G(f_{ex} )} \right| .
\]
If $\left| {D(\varphi _{ex} )} \right|<\varepsilon ^q$ then
\[
\left| {\frac{{H(g_{ex} )D(\varphi _{ex} )}} {{D^2 (\varphi _{ex} )
+ \delta _\varepsilon  }} - \frac{{H(g_{ex} )}} {{D(\varphi _{ex}
)}}} \right| \leqslant \chi \left( {B(\varphi _{ex},R_\varepsilon
,\varepsilon ^q)} \right).\left| {G(f_{ex} )} \right|\le \chi \left(
{B(\varphi _{ex},R_\varepsilon ,\varepsilon^q )} \right).\left\|
{f_{ex} } \right\|_{L^1 (\Omega )} ,
\]
where $B(\varphi _{ex} ,R_\varepsilon,\varepsilon^q)$ is as in Definition 1.
\\\h In summary, for all $(\alpha,n)\in R \times Z$ , we have
\[
\begin{gathered}
  \left| {G(f_\varepsilon  ) - G(f_{ex} )} \right| \leqslant 3C_3 \delta _\varepsilon ^{ - 1} \left| {H(g_\varepsilon   - g_{ex} )} \right| + 2\varepsilon \delta _\varepsilon ^{ - 1} \left| {H(g_{ex} )} \right|+\chi (B(\varphi _{ex} ,R_\varepsilon  ,\varepsilon ^q ))\left\| {f_{ex} } \right\|_{L^1 (\Omega )}  \hfill \\
\h\h\h\h\h\h+ \left[ {\chi ((\alpha ,n) \notin ( - R_\varepsilon  ,R_\varepsilon  )^2 ) + \chi ((\alpha ,n) \in ( - R_\varepsilon  ,R_\varepsilon  )^2 )\delta _\varepsilon  \varepsilon ^{ - 2q} } \right].\left| {G(f_{ex} )} \right|. \hfill \\
\end{gathered}
\]
To estimate $\sum\limits_{n =  - \infty }^\infty  {\int\limits_{ - \infty }^\infty  {\left| {G(f_\varepsilon  ) - G(f_{ex} )} \right|^2 d\alpha } }
$, we square the latter inequality to get
\[
\begin{gathered}
  \left| {G(f_\varepsilon  ) - G(f_{ex} )} \right|^2  \leqslant 36C_3^2 \delta _\varepsilon ^{ - 2} \left| {H(g_\varepsilon   - g_{ex} )} \right|^2  + 16\varepsilon ^2 \delta _\varepsilon ^{ - 2} \left| {H(g_{ex} )} \right|^2\hfill\\
\h\h\h\h\h\h   + 4\chi (B(\varphi _{ex} ,R_\varepsilon  ,\varepsilon ^q ))\left\| {f_{ex} } \right\|_{L^1 (\Omega )}^2\hfill \\
\h\h\h\h\h\h  + \left[ {4\chi ((\alpha ,n) \notin ( - R_\varepsilon  ,R_\varepsilon  )^2 ) + 4\delta _\varepsilon ^2 \varepsilon ^{ - 4q} } \right].\left| {G(f_{ex} )} \right|^2.  \hfill \\
\end{gathered}
\]
We will consider each term of the right-hand side. Since
\[
\left| {H(g_{ex} )} \right|^2  \leqslant \left( {\left| {G(g_{1ex} )} \right| + \left| {G(g_{0ex} )} \right|} \right)^2  \leqslant 2\left| {G(g_{1ex} )} \right|^2  + 2\left| {G(g_{0ex} )} \right|^2,
\]
we get
\[
\sum\limits_{n =  - \infty }^\infty  {\int\limits_{ - \infty }^\infty  {\left| {H(g_{ex} )} \right|^2 d\alpha } }  \leqslant 2\pi \left( {\left\| {g_{1ex} } \right\|_{L^2 (\Omega )}^2  + \left\| {g_{0ex} } \right\|_{L^2 (\Omega )}^2 } \right) \leqslant 4\pi C_3^2 .
\]
Similarly,
\[
\sum\limits_{n =  - \infty }^\infty  {\int\limits_{ - \infty }^\infty  {\left| {H(g_\varepsilon   - g_{ex} )} \right|^2 d\alpha } }  \leqslant 2\pi \left( {\left\| {g_{1\varepsilon }  - g_{1ex} } \right\|_{L^2 (\Omega )}^2  + \left\| {g_{0\varepsilon }  - g_{0ex} } \right\|_{L^2 (\Omega )}^2 } \right) \leqslant 4\pi \varepsilon ^2 .
\]
Furthermore,
\[
\begin{gathered}
  \sum\limits_{n =  - \infty }^\infty  {\int\limits_{ - \infty }^\infty  {\chi (B(\varphi _{ex} ,R_\varepsilon  ,\varepsilon ^q ))\left\| {f_{ex} } \right\|_{L^1 (\Omega )}^2 d\alpha } }  = m(B(\varphi _{ex} ,R_\varepsilon  ,\varepsilon ^q ))\left\| {f_{ex} } \right\|_{L^1 (\Omega )}^2,  \hfill \\
  \sum\limits_{n =  - \infty }^\infty  {\int\limits_{ - \infty }^\infty  {\chi ((\alpha ,n) \notin ( - R_\varepsilon  ,R_\varepsilon  )^2 )\left| {G(f_{ex} )} \right|^2 d\alpha } }  = \mu (f_{ex} ,R_\varepsilon  ), \hfill \\
  \sum\limits_{n =  - \infty }^\infty  {\int\limits_{ - \infty }^\infty  {\left| {G(f_{ex} )} \right|^2 d\alpha } }  = \pi \left\| {f_{ex} } \right\|_{L^2 (\Omega )}^2 . \hfill \\
\end{gathered}
\]
Thus
\[
\begin{gathered}
  \sum\limits_{n =  - \infty }^\infty  {\int\limits_{ - \infty }^\infty  {\left| {G(f_\varepsilon  ) - G(f_{ex} )} \right|^2 d\alpha } } \hfill\\
\leqslant 144\pi C_3^2 \varepsilon ^2 \delta _\varepsilon ^{ - 2}  + 64\pi C_3^2 \varepsilon ^2 \delta _\varepsilon ^{ - 2}  + 4m(B(\varphi _{ex} ,R_\varepsilon  ,\varepsilon ^q ))\left\| {f_{ex} } \right\|_{L^1 (\Omega )}^2\hfill\\
~~  + 4\mu (f_{ex} ,R_\varepsilon  ) + 4\pi \delta _\varepsilon ^2 \varepsilon ^{ - 4q} \left\| {f_{ex} } \right\|_{L^2 (\Omega )}^2 . \hfill \\
\end{gathered}
\]
Choosing $q=2/5$, $\delta_\varepsilon=\varepsilon^{9/10}$, we obtain that
\bq
\begin{gathered}
  \left\| {f_\varepsilon   - f_{ex} } \right\|_{L^2 (\Omega )}^2
  = \pi ^{ - 1} \sum\limits_{n =  - \infty }^\infty  {\int\limits_{ - \infty }^\infty  {\left| {G(f_\varepsilon  ) - G(f_{ex} )} \right|^2 d\alpha } }  \hfill \\
   \leqslant C_4 \varepsilon ^{1/5}  + 4\pi ^{ - 1} m(B(\varphi _{ex} ,R_\varepsilon  ,\varepsilon ^q ))\left\| {f_{ex} } \right\|_{L^1 (\Omega )}^2  + 4\pi ^{ - 1} \mu (f_{ex} ,R_\varepsilon  ), \hfill \\
\end{gathered} \label{xxss}
\eq
where $C_4  = 208 C_3^2  + 4\left\| {f_{ex} } \right\|_{L^2 (\Omega )}^2$.
\vspace{0.1in}\text{}\\
{\bf Step 3.} Proof of Theorem 2.
\\\h We construct $f_{1\varepsilon}$ as $f_{\varepsilon}$ corresponding
$$
q=2/5,\delta_\varepsilon=\varepsilon^{9/10},R_{\varepsilon}= R_{1\varepsilon}=(ln(\varepsilon^{-1}))^{\beta}.
$$
Applying $(i)$ of Lemma 2, for $\varepsilon>0$ small enough (depended on $\varphi_{ex}$ and $\beta$), we have $$m(B(\varphi _{ex} ,R_{1\varepsilon } ,\varepsilon ^q )) \leqslant R_{1\varepsilon }^{ - 1}.$$
Hence, from $(\ref{xxss})$, it follows that
\[
\left\| {f_{1\varepsilon}   - f_{ex} } \right\|_{L^2 (\Omega )}^2  \leqslant C_4 \varepsilon ^{1/5}  + 4\pi ^{ - 1} \left\| {f_{ex} } \right\|_{L^1 (\Omega )}^2 R_{1\varepsilon }^{ - 1}  + 4\pi ^{ - 1} \mu (f_{ex} ,R_{1\varepsilon } ).
\]
Since $\mathop {\lim }\limits_{\varepsilon  \to 0} R_{1\varepsilon }  =  + \infty$, Lemma 5 implies that $\mathop {\lim }\limits_{\varepsilon  \to 0} \mu (f_{ex} ,R_{1\varepsilon } ) = 0$. Therefore,
\[
\mathop {\lim }\limits_{\varepsilon  \to 0} \left\| {f_{1\varepsilon}   - f_{ex} } \right\|_{L^2 (\Omega )}  = 0.
\]
Moreover, if $f_{ex}\in H^1(\Omega)$ then Lemma 5 implies
\[
\mu (f_{ex} ,R_{1\varepsilon } ) \leqslant \left( {8R_{1\varepsilon }^{ - 1}  + 2\pi R_{1\varepsilon }^{ - 2} } \right)\left\| {f_{ex} } \right\|_{H^1 (\Omega )}^2 .
\]
So we obtain
\[
\begin{gathered}
  \left\| {f_{1\varepsilon}   - f_{ex} } \right\|_{L^2 (\Omega )}^2
\leqslant C_4 \varepsilon ^{1/5}  + 4\pi ^{ - 1} \left\| {f_{ex} } \right\|_{L^1 (\Omega )}^2 R_{1\varepsilon }^{ - 1}  + 4\pi ^{ - 1} \left( {8R_{1\varepsilon }^{ - 1}  + 2\pi R_{1\varepsilon }^{ - 2} } \right)\left\| {f_{ex} } \right\|_{H^1 (\Omega )}^2  \hfill \\
 \h\h\h\h~~\leqslant 4\pi ^{ - 1} \left( {1 + 9\left\| {f_{ex} } \right\|_{H^1 (\Omega )}^2 } \right)R_{1\varepsilon }^{ - 1}\hfill\\
 \h\h\h\h~~ = 4\pi ^{ - 1} \left( {1 + 9\left\| {f_{ex} } \right\|_{H^1 (\Omega )}^2 } \right)(ln(\varepsilon^{-1}))^{-\beta}.\hfill \\
\end{gathered}
\]
for $\varepsilon>0$ small enough (depended on $u_{ex},\varphi_{ex}$ and $\beta$).
\vspace{0.1in}\text{}\\
{\bf Step 4.} Proof of Theorem 3.
\\\h We construct $f_{2\varepsilon}$ as $f_{\varepsilon}$ corresponding
$$
q=2/5,\delta_\varepsilon=\varepsilon^{9/10},R_{\varepsilon}= R_{2\varepsilon}=\varepsilon^{1/6}.
$$
Applying $(ii)$ of Lemma 2 to $q_1=1/3$, there exists a constant $\gamma>0$ such that for $\varepsilon>0$ small enough (depended on $\varphi_{ex}$) we get
$$m(B(\varphi _{ex} ,R_{2\varepsilon } ,\varepsilon ^q )) \leqslant\varepsilon^{ \gamma}.$$
Hence, from $(\ref{xxss})$, it follows that
\[
\left\| {f_\varepsilon   - f_{ex} } \right\|_{L^2 (\Omega )}^2  \leqslant C_4 \varepsilon ^{1/5}  + 4\pi ^{ - 1} \varepsilon ^\gamma  \left\| {f_{ex} } \right\|_{L^1 (\Omega )}^2  + 4\pi ^{ - 1} \mu (f_{ex} ,R_{2\varepsilon } ).
\]
Since $\mathop {\lim }\limits_{\varepsilon  \to 0} R_{2\varepsilon }  =  + \infty$, we use Lemma 5 to get $\mathop {\lim }\limits_{\varepsilon  \to 0} \mu (f_{ex} ,R_{2\varepsilon } ) = 0$. It implies
\[
\mathop {\lim }\limits_{\varepsilon  \to 0} \left\| {f_{2\varepsilon}   - f_{ex} } \right\|_{L^2 (\Omega )}  = 0.
\]
Moreover, if $f\in H^1(\Omega)$ then Lemma 5 implies that
\[
\mu (f_{ex} ,R_{2\varepsilon } ) \leqslant \left( {8R_{2\varepsilon }^{ - 1}  + 2\pi R_{2\varepsilon }^{ - 2} } \right)\left\| {f_{ex} } \right\|_{H^1 (\Omega )}^2  = \left( {8\varepsilon ^{1/6}  + 2\pi \varepsilon ^{1/3} } \right)\left\| {f_{ex} } \right\|_{H^1 (\Omega )}^2 .
\]
Thus
\[
\begin{gathered}
  \left\| {f_\varepsilon   - f_{ex} } \right\|_{L^2 (\Omega )}^2  \leqslant C_4 \varepsilon ^{1/5}  + 4\pi ^{ - 1} \left\| {f_{ex} } \right\|_{L^1 (\Omega )}^2 \varepsilon ^\gamma   + 4\pi ^{ - 1} \left( {8\varepsilon ^{1/6}  + 2\pi \varepsilon ^{1/3} } \right)\left\| {f_{ex} } \right\|_{H^1 (\Omega )}^2  \hfill \\
  \h\h\h\h~ \leqslant 4\pi ^{ - 1} \left( {\left\| {f_{ex} } \right\|_{L^1 (\Omega )}^2 \varepsilon ^\gamma   + 9\left\| {f_{ex} } \right\|_{H^1 (\Omega )}^2 \varepsilon ^{1/6} } \right) .\hfill \\
\end{gathered}
\]
for $\varepsilon>0$ small enough (depended on $u_{ex}$ and $\varphi_{ex}$).
\end{proof}
\text{}\\
{\bf 3. A numerical experiment}
\\\\ We consider the exact data
\[
\begin{gathered}
  \varphi _{ex}  = \pi ^2 e^{ - \pi ^2 t} , \hfill \\
  g_{0ex}  = (\cos (\pi x) + 1)\cos (\pi y), \hfill \\
  g_{1ex}  = e^{ - \pi ^2 } (\cos (\pi x) + 1)\cos (\pi y). \hfill \\
\end{gathered}
\]
Then the corresponding exact solution of the system $(1)$ is
\[
\begin{gathered}
  u_{ex}  = e^{ - \pi ^2 t} (\cos (\pi x) + 1)\cos (\pi y), \hfill \\
  f_{ex}  = \cos (\pi x)\cos (\pi y) .\hfill \\
\end{gathered}
\]
For all $m=2,4,6,8,...$, we consider the disturbed data
\[
\begin{gathered}
  \varphi _m  = \varphi _{ex} , \hfill \\
  g_{0m}  = g_{0ex}  + m^{ - 1} (\cos (m\pi x) - 1)\cos (\pi y), \hfill \\
  g_{1ex}  = g_{1ex}  + m^{ - 1} (\cos (m\pi x) - 1)\cos (\pi y). \hfill \\
\end{gathered}
\]
Then the corresponding disturbed solution of the system $(1)$ is
\[
\begin{gathered}
  \widetilde u_m  = u_{ex}  + m^{ - 1} e^{ - \pi ^2 t} (\cos (m\pi x) - 1)\cos (\pi y), \hfill \\
  \widetilde f_m  = f_{ex}  + m\cos (m\pi x)\cos (\pi y) .\hfill \\
\end{gathered}
\]
We get
\[
\begin{gathered}
  \left\| {g_{0m}  - g_{0ex} } \right\|_{L^2 (\Omega )}  = \frac{{\sqrt 3 }}
{{2m}},  \hfill \\
 \left\| {g_{0m}  - g_{0ex} } \right\|_{L^2 (\Omega )}  = \frac{{e^{ - \pi ^2 } \sqrt 3 }}
{{2m}} ,\hfill \\
  \left\| {\widetilde f_m  - f_{ex} } \right\|_{L^2 (\Omega )}  = \frac{m}
{2}. \hfill \\
\end{gathered}
\]
It means that, when $m$ is large, a small error of data causes a large error of solution. Hence, the problem is ill-posed and a regularization is necessary.
\\\h We construct $f_{1m}$ and $f_{2m}$ as $f_{1\varepsilon}$ and $f_{2\varepsilon}$ in Theorem 2 corresponding
$$
\begin{gathered}
  \varepsilon  = m^{ - 1} ,\delta _\varepsilon   = \varepsilon ^{9/10}  = m^{ - 9/10} ,\hfill \\
  \varphi _\varepsilon   = \varphi _m ,g_{0\varepsilon }  = g_{0m} ,g_{1\varepsilon }  = g_{1m} , \hfill \\
R_{1\varepsilon }  =  R_{1m}=(ln(m))^{2/5},R_{2\varepsilon }  =  R_{2m}=m^{1/6}.  \hfill \\
\end{gathered}
$$
We have
\[
\begin{gathered}
  D(\varphi _m )(\alpha ,n) = \frac{{\pi ^2 e^{ - \pi ^2 } (1 - e^{ - \alpha ^2 } )}}
{{\alpha ^2 }}, \hfill \\
  H(g_m )(\alpha ,n) = G(g_{1m} )(\alpha ,n) - e^{ - (\alpha ^2  + n^2 \pi ^2 )} G(g_{0m} )(\alpha ,n) \hfill \\
   = \left\{ \begin{gathered}
  0\text{ if } n\ne \pm 1, \hfill \\
  \frac{{(m - 1)\pi ^2 e^{ - \pi ^2 } (1 - e^{ - \alpha ^2 } )\sin (\alpha )(\alpha ^2  + m\pi ^2 )}}
{{2\alpha (\alpha ^2  - \pi ^2 )(\alpha ^2  - m^2 \pi ^2 )}}\text{ if } n=\pm 1.\hfill \\
\end{gathered}  \right. \hfill \\
\end{gathered}
\]
Hence, for each $j\in\{1,2\}$, we get
\[
\begin{gathered}
  f_{jm} (x,y) = \frac{1}
{\pi }.\sum\limits_{ - R_{jm}  < n < R_{jm} } {\left( {\int\limits_{ - R_{jm} }^{R_{jm} } {\frac{{H(g_m ).D(\varphi _m )}}
{{D^2 (\varphi _m ) + m^{ - 9/10} }}\cos (\alpha x)d\alpha } } \right)\cos (n\pi y)}  \hfill \\
   = \frac{2}
{\pi }.\left( {\int\limits_{ - R_{jm} }^{R_{jm} } {\frac{{H(g_m )(\alpha ,1).D(\varphi _m )(\alpha ,1)}}
{{D^2 (\varphi _m )(\alpha ,1) + m^{ - 9/10} }}\cos (\alpha x)d\alpha } } \right)\cos (\pi y) \hfill \\
   = \frac{2}
{\pi }.\left( {\int\limits_{ - R_{jm} }^{R_{jm} } {\frac{{(m - 1)\pi ^4 \alpha \sin (\alpha )(\alpha ^2  + m\pi ^2 )}}
{{2(\alpha ^2  - \pi ^2 )(\alpha ^2  - m^2 \pi ^2 )\left( {\pi ^4  + m^{ - 9/10} \alpha ^4 e^{2\pi ^2 } \left( {1 - e^{ - \alpha ^2 } } \right)^{ - 2} } \right)}}\cos (\alpha x)d\alpha } } \right)\cos (\pi y). \hfill \\
\end{gathered}
\]
\h The error estimates between $f_{jm}$ and $f_{ex}$ is given by the following table.
\newpage
\begin{center}
\begin{tabular}{|l|l|l|}
\hline $\varepsilon=m^{-1}$ & $\left\| {f_{1m}  - f_{ex} } \right\|_{L^2 (\Omega )}^2$ & $\left\| {f_{2m}  - f_{ex} } \right\|_{L^2 (\Omega )}^2$\\ \hline
$10^{-2}$ & 0.2499999413 & 0.2499999226 \\
$10^{-6}$ & 0.2495687285 & 0.2494516526 \\
$10^{-12}$ & 0.1181388651 & $1.240475046 \times 10^{-2}$\\
$10^{-15}$ & $9.038948894\times 10^{-2}$ & $5.031837329 \times 10^{-4}$ \\
$10^{-20}$ & $5.881288742 \times 10^{-2}$ & $7.385640149 \times 10^{-5}$ \\
$10^{-30}$ & $3.361360762 \times 10^{-2}$ & $1.591595841  \times 10^{-6}$ \\
 \hline
\end{tabular}
\vspace{0.1in}\text{ }\\Table 1. Error estimates between $f_{jm}$ and $f_{ex}$
\end{center}
Here are some figures of the exact solution and the regularized
solution with $m=10^{15}$.\ \par \vspace{0.2in}
 \centerline{
\includegraphics[width=4in]{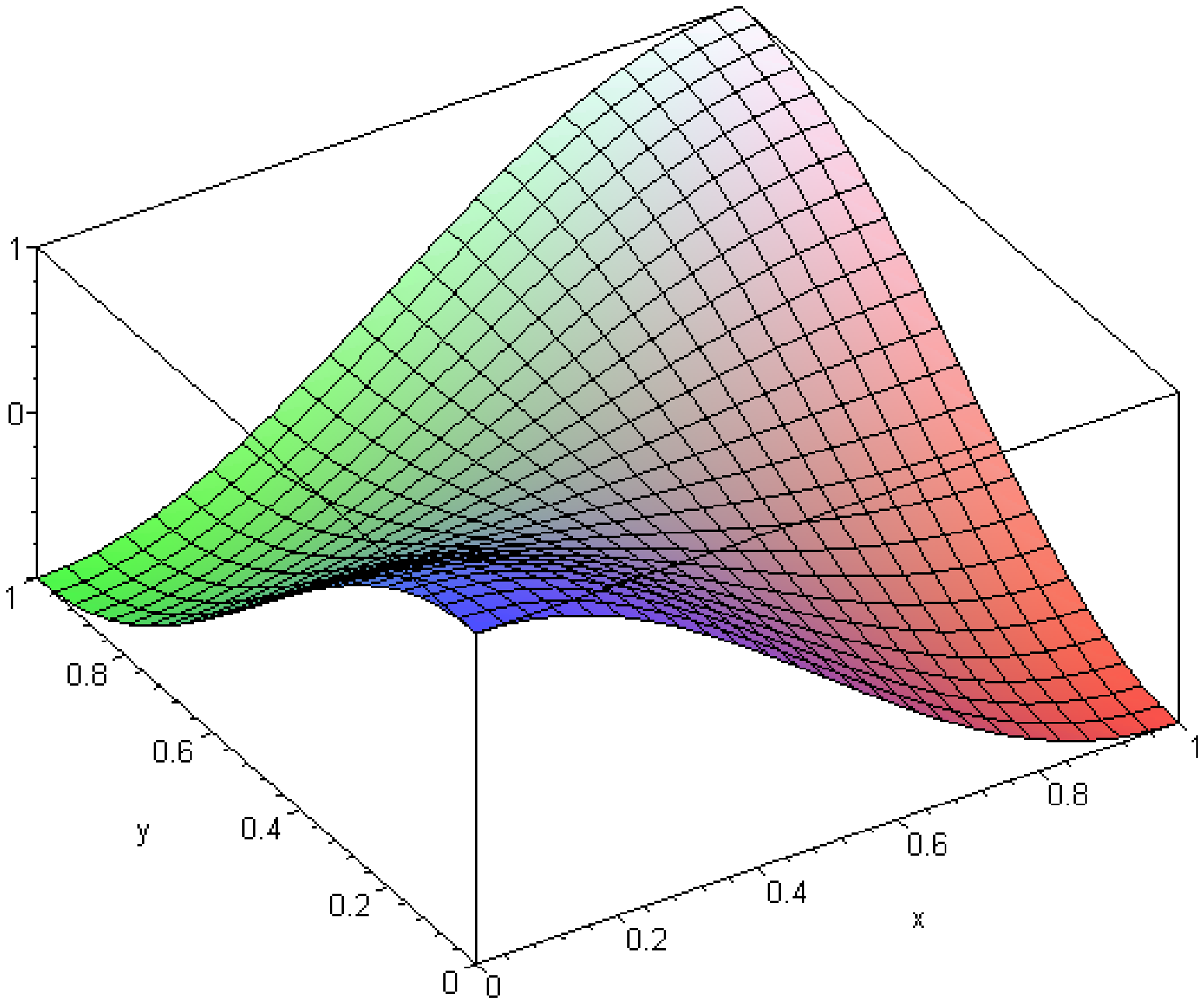}}
\center{Figure 1. The exact solution}\endcenter
\ \par \centerline{
\includegraphics[width=4in]{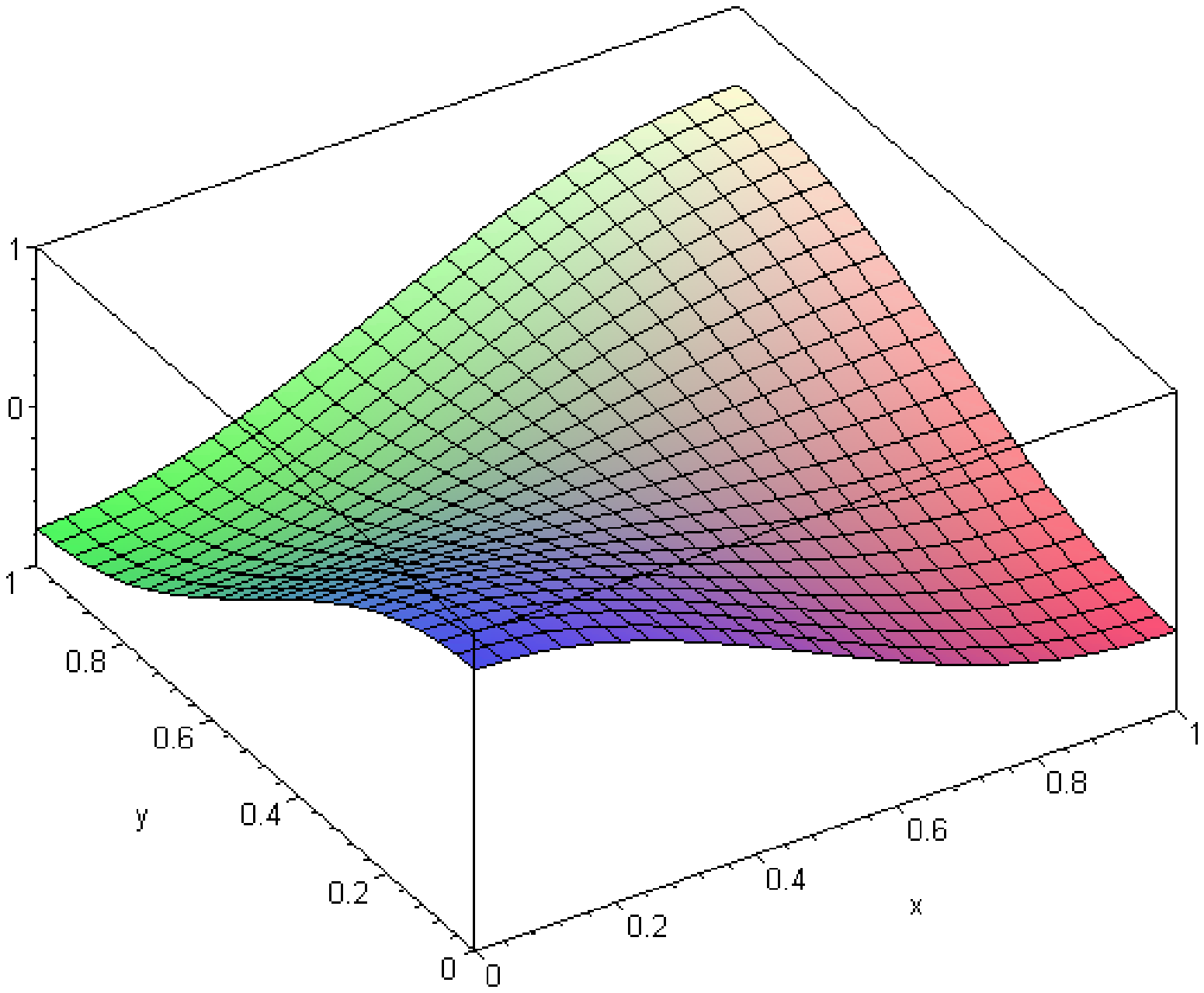}}
\center{Figure 2. The regularized solution $f_{1m}$}\endcenter
\newpage \centerline{
\includegraphics[width=4in]{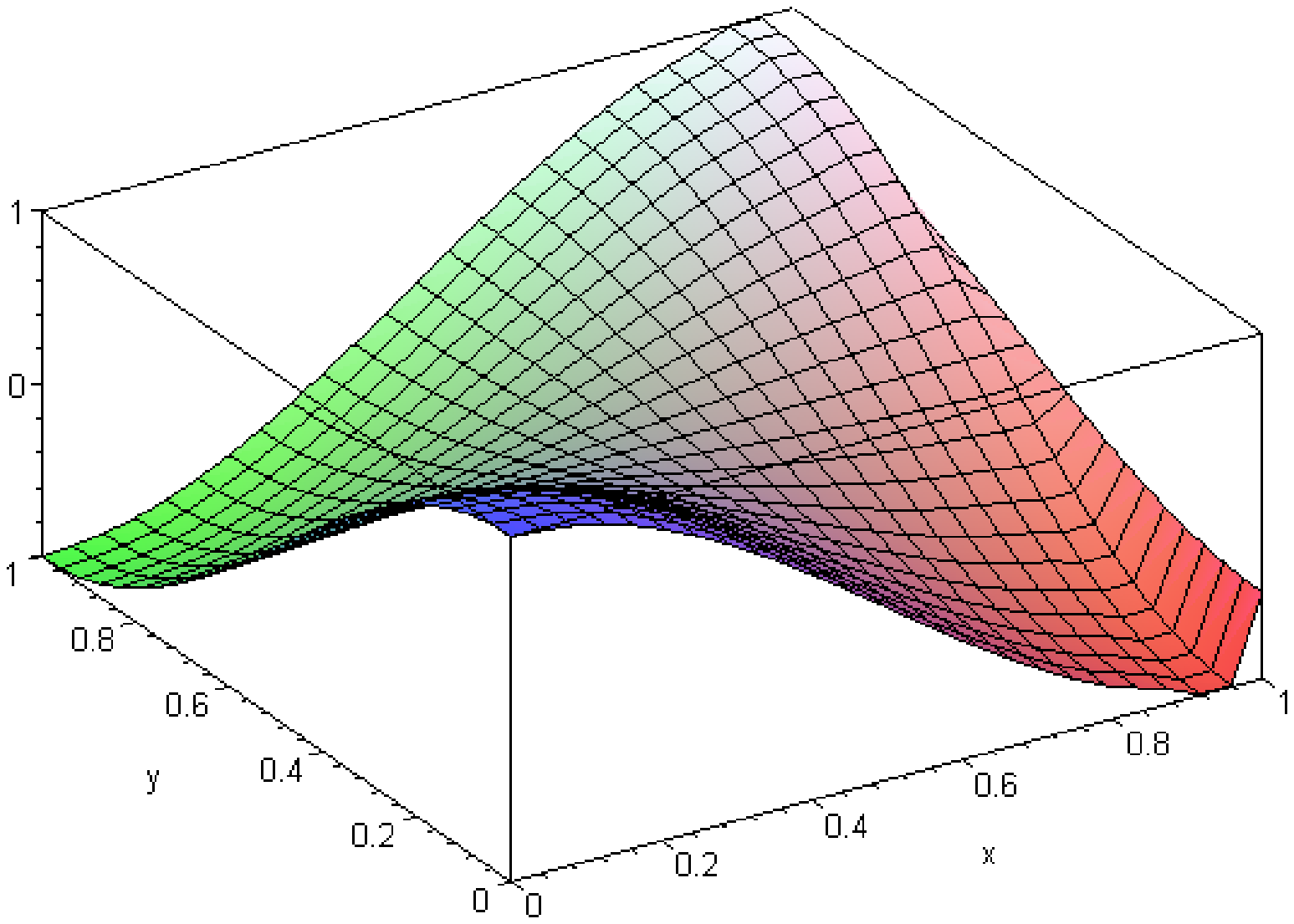}}
\center{Figure 3. The regularized solution $f_{2m}$ }\endcenter
\begin{remark} $f_{jm}$ approximates $f_{ex}$ well on $\Omega$ except a neighborhood of the  line $y=1$. The fact is understandable because the value of $u(x,1,t)$ is not given in the system $(1)$.
\end{remark}

\end{document}